\documentclass[12pt]{amsart}
\usepackage{amsmath}
\usepackage{amssymb, amsfonts, amsthm}
\usepackage{mathtools}
\usepackage{amssymb}
\usepackage{ifthen}
\usepackage{graphicx}
\usepackage{float}
\usepackage{caption}
\usepackage{subcaption}
\usepackage{cite}
\usepackage{mathrsfs}
\usepackage{amsfonts}
\usepackage{amscd}
\usepackage{amsxtra}
\usepackage{color}
\usepackage{bm}
\numberwithin{equation}{section}

\newtheorem{theorem}{Theorem}[section]
\newtheorem{lemma}[theorem]{Lemma}

\theoremstyle{definition}
\newtheorem{definition}[equation]{Definition}
\newtheorem{example}[theorem]{Example}

\newtheorem{remark}[equation]{Remark}

\newcounter{alphabet}
\newcounter{tmp}
\newenvironment{Thm}[1][]{\refstepcounter{alphabet}%
	\bigskip%
	\noindent%
	{\bf Theorem \Alph{alphabet}}%
	\ifthenelse{\equal{#1}{}}{}{ (#1)}%
	{\bf .} \itshape}{\vskip 8pt}

\makeatletter

\newcommand{\RRef}[1]{\@ifundefined{r@#1}{}{\setcounter{tmp}{\ref{#1}}\Alph{tmp}}}

\newcounter{minutes}
\divide\time by 60
\newcounter{hours}
\multiply\time by 60 \addtocounter{minutes}{-\time}

\begin{document}
	\title[On the Boundary Schwarz lemma and the rigidity theorem for certain mappings]
	{On the Boundary Schwarz lemma and the rigidity theorem for certain mappings}
	
	\thanks{
		File:~\jobname .tex,
		printed: \number\day-\number\month-\number\year,
		\thehours.\ifnum\theminutes<10{0}\fi\theminutes}
	
	\author[Shankey Kumar]{Shankey Kumar}
	\address{Shankey Kumar, Department of Mathematics, Indian Institute of Technology Madras, Chennai, 600036, India.}
	\email{shankeygarg93@gmail.com}
	
	\author[Saminathan Ponnusamy]{Saminathan Ponnusamy}
	\address{Saminathan Ponnusamy, Department of Mathematics, Indian Institute of Technology Madras, Chennai, 600036, India.}
	\email{samy@iitm.ac.in}
	
	
	\begin{abstract}
In this article, we characterize the holomorphic mappings from $B_{\ell_p^n}\times\mathbb{D}^{m}$ into $\mathbb{D}^{m}$ for $p\in \{2,\infty\}$. In addition, we give a simple proof for the boundary Schwarz lemma for vector valued holomorphic functions, which also extends the existing result. Also, we obtain the boundary Schwarz lemma for pluriharmonic self-mappings of the unit ball $B_{\ell_p^n}$, $p \in [2,\infty]$. Furthermore, we establish the boundary rigidity theorem for holomorphic self-mappings of $B_{\ell_p^n}$, $p \in (1,\infty)$.
	\end{abstract}
	
	\maketitle
	\pagestyle{myheadings}
	\markboth{S. Kumar and S. Ponnusamy}{THE BOUNDARY SCHWARZ LEMMA AND THE RIGIDITY THEOREM}
	
\section{Introduction}

The Schwarz lemma, together with its many generalizations, is fundamental in solving many extremal problems in geometric function theory and hyperbolic geometry. Because of its importance, the Schwarz lemma has been studied by many researchers (see also  
the recent monographs by Defant et al. \cite{ADM2019}  and Garcia \cite{GarMasRoss-2018} in which many challenging research problems have been discussed). In the recent years, special attention has been given to the behavior of the Schwarz lemma near the boundary of domains. Here is the one-variable version:

\begin{Thm}\label{BSL}{\em\cite{G81,G98}}
Let $f$ be holomorphic in the unit disk  $\mathbb{D}=\{z \in \mathbb{C}:\, |z| < 1\}$. If $f$ is holomorphic at $z=1$ with $f(0)=0$ and $f(1)=1$, then
	$f^{\prime}(1) \geq 1$, with equality $f^{\prime}(1)=1$ holds if and only if $f(z)=z$.
\end{Thm}

The boundary Schwarz lemma also has an important role in classical complex analysis and its applications.
In addition, the boundary Schwarz lemma is a key tool in studying bounded domains in $\mathbb{C}^{n}$. See \cite{BK94, GHK20, 
H95,  Z22}. In the case of the unit ball $B_{\ell_2^n}$ (definition is given below), Liu et al. \cite{LW15} (see Theorem~B) 
proved the  boundary Schwarz lemma for holomorphic self-maps of the unit ball.


More recently, a boundary rigidity theorem has been established for holomorphic mappings with values in finite-dimensional bounded symmetric domains (see \cite{HK21}). In addition, in \cite{WZ24}, the authors considered holomorphic mappings defined on the unit ball $B_{\ell_p^n}$ for $p \in [2, \infty]$, to study rigidity theorem.

The main purpose of this paper is to extend the boundary Schwarz lemma for holomorphic and pluriharmonic mappings in certain Reinhardt domains.
We establish (see Theorems~\ref{thm3.1}) a boundary Schwarz lemma for holomorphic self-mappings of $B_{\ell_p^n}$ with $p \in [2,\infty)$.
We further extend the work of Kalaj \cite{K25} to Banach space-valued holomorphic mappings and provide a simple proof of Theorem E with sharp estimate.
In Theorem \ref{thmph}, we generalize a result of Hamada \cite[Theorem 1.2]{H17} to Reinhardt domains in $\mathbb{C}^n$.
Moreover, Wang and Zhang \cite{WZ24} established a rigidity theorem (see Theorem K) for holomorphic functions on $B_{\ell_p^n}$, $p\in[2,\infty]$.
In Theorem \ref{rigi1}, we address the remaining cases and obtain the rigidity theorem for $p \in (1,\infty)$.

Several of the well-known results of different authors will be recalled in Section \ref{background}.
In order to state and prove our main results, we need to fix up some notations and recall some known results.

\subsection{Notations.}\label{notation}
The spaces $\mathbb{C}^n$ is the $n$-dimensional complex Hilbert space with the standard inner product
$ \langle z, w \rangle = \sum_{j=1}^{n} z_j \overline{w_j},
$
where $z=(z_1,\dots,z_n)$ and $w=(w_1,\dots,w_n) \in \mathbb{C}^n$. The unit polydisk in $\mathbb{C}^n$ is
$$\mathbb{D}^n = \{ z \in \mathbb{C}^n : \|z\|_\infty < 1 \},$$ where $\|z\|_\infty = \max \{ |z_j| : 1 \leq j \leq n \}$.
For $p \geq 1$, we define
$$
B_{\ell_p^n} = \Big\{ (z_1,\dots,z_n) \in \mathbb{C}^n :
\|z\|_p = (|z_1|^p + \cdots + |z_n|^p)^{1/p} < 1 \Big\}.
$$
For $p=2$, $B_{\ell_p^n}$ is the Euclidean unit ball.
Similarly, let $B_p^{2n}$ be the unit ball in $\mathbb{R}^{2n}$.
Each point $z = x + i y \in \mathbb{C}^n$ corresponds to
$
z' = (x,y)^T \in \mathbb{R}^{2n},
$
where $T$ denotes for transpose. It is customary to treat $z \in \mathbb{C}^n$ either as a column vector
$ z = (z_1, z_2, \dots, z_n)^T,
$
or as $n \times 1$ matrix
$$
\begin{pmatrix}
	 z_1  \\
	\vdots \\
 z_n
\end{pmatrix}.
$$

Let $\Omega_1$ and $\Omega_2$ denote two open sets in $\mathbb{C}^n$.
Denote by $H(\Omega_1,\Omega_2)$ the set of all holomorphic maps from $\Omega_1$ into $\Omega_2$.
If $\Omega_1 = \Omega_2$, we simply write $H(\Omega_1)$, instead of $H(\Omega_1, \Omega_1)$.
For $f \in H(\Omega_1, \Omega_2)$, we write
$$
f = (f_1, \dots, f_n)^T,
$$
where each $f_j : \Omega_1 \to \mathbb{C}$ is holomorphic.
The derivative of $f$ at $a \in \Omega_1$ is the complex Jacobian matrix
$$
J_f(a) = \left( \frac{\partial f_i}{\partial z_j}(a) \right)_{n \times n}
= \begin{pmatrix}
	\frac{\partial f_1}{\partial z_1}(a) & \cdots & \frac{\partial f_1}{\partial z_n}(a) \\
	\vdots & \ddots & \vdots \\
	\frac{\partial f_n}{\partial z_1}(a) & \cdots & \frac{\partial f_n}{\partial z_n}(a)
\end{pmatrix}.
$$
Thus $J_f(a)$ is a linear map from $\mathbb{C}^n$ to $\mathbb{C}^n$. We write $\|J_f(a)\|$ for the operator norm of $J_f(a)$. For $f:\Omega_1 \to \Omega_2$, where $\Omega_1\in \mathbb{C}^n$ and $\Omega_2\in \mathbb{C}^m$, denote $J_f(z')$ by the $2n\times 2m$ Jacobian matrix of $f$ at $z$ in terms of real coordinates.

A $C^2$ map $f : B_{\ell_p^n} \to \mathbb{C}^n$ is called {\em pluriharmonic} if the restriction of each component $f_j$ to any complex line is harmonic.

Finally, if $\Omega \subset \mathbb{R}^n$ and if $f : \Omega \to \mathbb{R}^n$ is $C^1$, then the Jacobian of $f$ at $x \in \Omega$ is the $n \times n$ matrix $J_f(x)$ with real entries.

\subsection{Background}\label{background}
%
%
%

The following version of the boundary Schwarz lemma for holomorphic self-maps of the unit ball  can be seen as a higher-dimensional version of Theorem A. 

\begin{Thm}\label{LW} {\em\cite{LW15}}
	Let $f \in H(B_{\ell_2^n})$ be holomorphic at some boundary point $z_{0} \in \partial B_{\ell_2^n}$ with $f(z_{0})=z_{0}$. Then the following hold:
	\begin{enumerate}
		\item [(1)] There exists a real number $\lambda$ such that
		$
		\overline{J_{f}(z_{0})}^{T} z_{0} = \lambda z_{0}.
		$
		\item [(2)] We have
		$$
		\lambda = \overline{z_{0}}^{T} J_{f}(z_{0}) z_{0} \geq
		\frac{|1-\overline{f(0)}^{T}z_{0}|^{2}}{1-\|f(0)\|_{2}^{2}} > 0,
		$$
		where $\overline{J_{f}(z_{0})}^{T}$ is the transpose of the conjugate matrix of $J_{f}(z_{0})$.
		In particular, if $f(0)=0$, then $\lambda \geq 1$.
		\item [(3)] The Jacobian determinant satisfies
		$$
		\left|\det J_{f}(z_{0})\right| \leq \lambda^{\frac{n+1}{2}}.
		$$
	\end{enumerate}
\end{Thm}


In \cite{H17}, Hamada improved a result of Liu et al. \cite{LDP16} and established the following boundary Schwarz lemma for pluriharmonic mappings.

\begin{Thm}\label{Hamada} {\em\cite{H17}}
	Let $f: B_{\ell_2^n} \rightarrow B_{\ell_2^N}$ be a pluriharmonic mapping for $n, N \geq 1$. Assume that $f$ is $C^{1}$ at $z_{0} \in \partial B_{\ell_2^n}$ and $f\left(z_{0}\right)=w_{0} \in \partial B_{\ell_2^N}$.
	Then
	$$
	\left(J_{f}\left(z_{0}^{\prime}\right)z_{0}^{\prime}\right)^{T} w_{0}^{\prime} \geq \frac{1-\left(f(0)^{\prime}\right)^{T} w_{0}^{\prime}}{2} \geq \frac{1-\|f(0)\|_2}{2}>0.
	$$
\end{Thm}

Zhu \cite{Z18} obtained the following version of the boundary Schwarz lemma for holomorphic self-mappings of the unit disk $\mathbb{D}$.

\begin{Thm}\label{Zhu} {\em \cite{Z18}}
Let $f\in H(\mathbb{D})$ and $f$ be holomorphic at $z=1$ with $f(1)=1$. Then
\begin{equation}\label{Zeq}
	f^{\prime}(1) \geq \frac{2|1-f(0)|^{2}}{1-|f(0)|^{2}+\left|f^{\prime}(0)\right|}.
\end{equation}
The above inequality is sharp with the extremal function
\begin{equation*}
	\varphi(z)=\frac{\beta A(z)+f(0)}{1+\beta \overline{f(0)} A(z)},
\end{equation*}
where
$$
\beta=\frac{1-f(0)}{1-\bar{f}(0)} \in \partial\mathbb{D}
~\mbox{ and }~
A(z)=z \frac{\left(1-|f(0)|^{2}\right) z+\left|f^{\prime}(0)\right|}{\left(1-|f(0)|^{2}\right)+\left|f^{\prime}(0)\right| z}.
$$
\end{Thm}
Recently, in \cite{K25}, Kalaj extended the above result to mappings $f\in H(\mathbb{D}, B_{\ell_2^n})$.

\begin{Thm}\label{Kalaj} {\em\cite{K25}}
	Let $f\in H(\mathbb{D}, B_{\ell_2^n})$ and $f(1)\in \partial B_{\ell_2^n}$. If $f'(1)$ exists, then we have
	\begin{equation*}
		\|f'(1)\|_2 \geq \frac{2(1-\|f(0)\|_2)^2}{1-\|f(0)\|_2^2+\|f'(0)\|_2}.
	\end{equation*}
\end{Thm}

A very interesting part of the boundary Schwarz lemma (Theorem A) is the uniqueness statement, also known as the rigidity property of holomorphic functions. The first major rigidity theorem of this type is due to Burns and Krantz \cite{BK94}, which gives a boundary version of a rigidity result for holomorphic maps:

\begin{Thm} {\em \cite[Theorem 3.1]{BK94}}
	If $f \in H(B_{\ell_2^n})$ satisfies
	$$
	f(z)=\mathbf{1}+(z-\mathbf{1})+o(\|z-\mathbf{1}\|^{3})
	\quad \text{as } z \to \mathbf{1},
	$$
	where $\mathbf{1}=(1,0,\dots,0)$, then $f(z)=z$.
\end{Thm}

Zimmer \cite{Z22} recently proved a rigidity theorem for bounded convex domains with $C^{2}$ boundaries, using the Kobayashi metric:

\begin{Thm}{\em \cite[Theorem 1.5]{Z22}}
	Let $\Omega \subset \mathbb{C}^{n}$ be a bounded convex domain with $C^{2}$ boundary and $\xi_{0} \in \partial \Omega$. If $f \in H(\Omega)$ satisfies
	$$
	f(z)=z+o(\|z-\xi_{0}\|_{2}^{4}) \quad \text{as } z \to \xi_{0},
	$$
	then $f(z)=z$.
\end{Thm}

For strongly convex domains, Huang \cite{H95} obtained the following:

\begin{Thm} {\em\cite[Corollary 2.7]{H95}}
	Let $\Omega \subset \mathbb{C}^{n}$ be strongly convex with $n>1$. Suppose $p \in \partial \Omega$ and $z_{0} \in \Omega$. If $f \in H(\Omega)$ with $f(z_{0})=z_{0}$ and
	$$
	f(z)=z+o(\|z-p\|_{2}^{2}) \quad \text{as } z \to p,
	$$
	then $f(z)=z$.
\end{Thm}

As an application of Theorem B, Tang et al. \cite{TLZ17} proved a new rigidity result on the unit ball.

\begin{Thm} {\em\cite{TLZ17}}
	Let $f \in H(B_{\ell_2^n})$ with $f(0)=0$. Suppose there exist linearly independent vectors $\alpha_{1},\dots,\alpha_{n} \in \partial B_{\ell_2^n}$ such that $f$ is holomorphic at $z=\alpha_{k}$ and $f(\alpha_{k})=\alpha_{k}$ for $k=1,\dots,n$. Then
	$$
	\overline{\alpha_{k}}^{T} J_{f}(\alpha_{k}) \alpha_{k}=1
	\quad (k=1,\dots,n)
	$$
	if and only if $f(z)=z$.
\end{Thm}

For the unit polydisk $\mathbb{D}^{n}$, Tang and Liu \cite{LT20} proved a similar rigidity theorem.

\begin{Thm}{\em\cite{LT20}}
	Let $f \in H(\mathbb{D}^{n})$ with $f(0)=0$. Suppose there exist linearly independent vectors $\alpha_{1},\dots,\alpha_{n} \in (\partial \mathbb{D})^{n}$ such that $f$ is holomorphic at $z=\alpha_{k}$ and $f(\alpha_{k})=\alpha_{k}$ for $k=1,\dots,n$. Then
	$$
	\overline{\alpha_{k}}^{T} J_{f}(\alpha_{k}) \alpha_{k}=n
	\quad (k=1,\dots,n)
	$$
	if and only if $f(z)=z$.
\end{Thm}
Wang and Zhang \cite{WZ24} established the following rigidity result for holomorphic functions on $B_{\ell_p^n}$, which extends the results from the Euclidean unit ball (corresponding to $p=2$) and the unit polydisk (corresponding to $p=\infty$).

\begin{Thm}\label{rigiWZ} {\em\cite{WZ24}}
	Suppose that $p \in[2, \infty]$ and that $f\in H(B_{\ell_p^n})$ with $f(0)=0$. Let
	$$
	v_{z}= \begin{cases}\left(\left|z_{1}\right|^{p-2} z_{1}, \dots,\left|z_{n}\right|^{p-2} z_{n}\right)^{T} & \text { if } p \geq 2, z \in \partial B_{\ell_p^n}, \\ \cfrac{z}{n} & \text { if } p=\infty, z \in(\partial \mathbb{D})^{n}.\end{cases}
	$$
	
	If there exist linearly independent vectors $\alpha_{1}, \alpha_{2}, \dots, \alpha_{n} \in \partial B_{\ell_p^n}$ such that $f$ is holomorphic at $z=\alpha_{k}$ and $f\left(\alpha_{k}\right)=\alpha_{k}(k=1, \dots, n)$, then the $n$ equations
	$$
	{\overline{v_{\alpha_{k}}}}^{T} J_{f}\left(\alpha_{k}\right) \alpha_{k}=1 \quad(k=1, \dots, n)
	$$
	hold if and only if $f(z)=z$.
\end{Thm}

%
\section{Preliminaries and main results}
In this section, we recall several standard definitions and lemmas that will be used in the proofs of our main results. We also present the statements of the main theorems.

\subsection{Schwarz Lemma} The following result characterizes  the holomorphic maps from $B_{\ell_p^n}\times\mathbb{D}^{m}$ into $\mathbb{D}^{m}$ for $p\in \{2,\infty\}$:
\begin{theorem}\label{SL}
Let $p\in \{2,\infty\}$. If $f\in H(B_{\ell_p^n}\times\mathbb{D}^{m}, \mathbb{D}^{m})$, $\phi$ is an automorphism of $\mathbb{D}^{m}$, and there exists a $z_{0} \in B_{\ell_p^n}$ such that $f\left(z_{0}, w\right)=\phi(w)$ for all $w$, then $f(z, w)=\phi(w)$ for all $(z, w)$.
\end{theorem}


\begin{remark}
	The case $m=1$ and $p=\infty$ of Theorem \ref{SL} was established by Knese \cite{K11}.
\end{remark}

\subsection{Carath\'{e}odory Metric}
Let $\Omega$ be a domain in $\mathbb{C}^{n}$.
The Carath\'{e}odory metric $\mathcal{C}_{\Omega}: \Omega \times \mathbb{C}^{n} \to \mathbb{R}^{+}$ is defined by
$$\mathcal{C}_{\Omega}(z, \xi) = \sup_{ f \in \mathcal{H}(\Omega, \mathbb{D}), \ f(z)=0}
\left| \sum_{j=1}^{n} \frac{\partial f}{\partial z_{j}}(z)\,\xi_{j} \right|  .
$$
It is easy to see that the metric is homogeneous, i.e.,
$$
\mathcal{C}_{\Omega}(z, t\xi) = |t|\,\mathcal{C}_{\Omega}(z, \xi) \quad \text{for } t \in \mathbb{C}.
$$
For all $z, w \in \Omega$, the Carath\'{e}odory distance is defined by
$$d_{\Omega}^{c}(z, w) = \sup_{f \in \mathcal{H}(\Omega, \mathbb{D}),~f(z)=0}   
\omega (f(z),f(w)) ,
$$
where $\omega (a,b)$ is the hyperbolic distance on $\mathbb{D}$, namely,
$$ \omega (a,b):= d_{\mathbb{D}}^{c}(a, b) =
\frac{1}{2} \log \frac{|1-\bar{a} b| + |a-b|}{|1-\bar{a} b| - |a-b|}.
$$

For a detailed overview of this topic, see Krantz \cite{K01} and Gong \cite{G98}.
%
The following lemma contained in \cite[Proposition 2.2.1]{JP93} gives a Schwarz-type result for holomorphic self-maps of the unit ball $B_{\ell_p^n}$ for $p \in [1,\infty]$ (see also \cite{WZ24}).
\begin{lemma} \label{lem2.3}
For any $p\in [1,\infty]$ and $z \in B_{\ell_p^n}$, we have that
$$
\mathcal{C}_{B_{\ell_p^n}}(0, z)=\|z\|_{p} ~\mbox{ and }~
d_{B_{\ell_p^n}}^{c}(0, z)=\frac{1}{2} \log \frac{1+\|z\|_{p}}{1-\|z\|_{p}}.
$$
Furthermore, if $f\in H(B_{\ell_p^n})$ with $f(0)=0$, then
$$
\|f(z)\|_{p} \leq\|z\|_{p} ~\mbox{ and }~
\left\|f^{\prime}(0)\right\|=\max _{\|\xi\|_{p} \leq 1}\left\|f^{\prime}(0) \xi\right\|_{p}=\max _{\|\xi\|_{p}=1}\left\|f^{\prime}(0) \xi\right\|_{p} \leq 1.
$$
\end{lemma}

\subsection{Tangent Space}
Let $p \in[2, \infty)$ and let $\rho(z)=\|z\|_{p}^{p}-1=\left|z_{1}\right|^{p}+\cdots+\left|z_{n}\right|^{p}-1$. Then the gradient of $\rho(z)$ is
$$
\nabla \rho(z)=2\left(\frac{\partial \rho}{\partial \overline{z_{1}}}(z), \dots, \frac{\partial \rho}{\partial \overline{z_{n}}}(z)\right)^{T}=p\left(\left|z_{1}\right|^{p-2} z_{1}, \dots,\left|z_{n}\right|^{p-2} z_{n}\right)^{T}
$$
which is an outer normal vector at the boundary point $z_{0}$. Moreover, it is easy to see that $\rho(z)$ is a $C^{2}$ defining function of the domain $B_{\ell_p^n}$ for $p \in[2, \infty)$.

\begin{definition}
For any fixed $z \in \partial B_{\ell_p^n}$, the tangent space $\mathcal{T}_{z}\left(\partial B_{\ell_p^n}\right)$ to $\partial B_{\ell_p^n}$ at $z$ is defined by
$$
\mathcal{T}_{z}\left(\partial B_{\ell_p^n}\right)=\left\{\alpha \in \mathbb{C}^{n}: {\rm Re}\left\langle\alpha, \nabla \rho\left(z\right)\right\rangle={\rm Re}\left[{\overline{\nabla \rho\left(z\right)}}^{T} \alpha\right]=0\right\}.
$$
The complex tangent space $\mathcal{T}_{z}^{1,0}\left(\partial B_{\ell_p^n}\right)$ to $\partial B_{\ell_p^n}$ at $z \in \partial B_{\ell_p^n}$ is defined by
$$
\mathcal{T}_{z}^{1,0}\left(\partial B_{\ell_p^n}\right)=\left\{\alpha \in \mathbb{C}^{n}:\left\langle\alpha, \nabla \rho\left(z\right)\right\rangle={\overline{\nabla \rho\left(z\right)}}^{T} \alpha=0\right\}.
$$
\end{definition}

\subsection{ Boundary Schwarz Lemma}
 The following result generalizes Theorem A to  holomorphic mappings from the space $H(B_{\ell_p^n})$, $p \in[2, \infty)$, using the similar approach as  in \cite[Theorem 3.1]{WZ24}:

\begin{theorem} \label{thm3.1} For $p \in[2, \infty)$, let $f\in H(B_{\ell_p^n})$ with $f(0)=0$. If $f$ is holomorphic at $z_{0} \in \partial B_{\ell_p^n}$ and $f\left(z_{0}\right)=w_{0} \in \partial B_{\ell_p^n}$, then the following two statements hold:
	\begin{itemize}
		\item[(1)] $J_{f}\left(z_{0}\right) \mathcal{T}_{z_{0}}\left(\partial B_{\ell_p^n}\right) \subset \mathcal{T}_{w_{0}}\left(\partial B_{\ell_p^n}\right)$ and $J_{f}\left(z_{0}\right) \mathcal{T}_{z_{0}}^{1,0}\left(\partial B_{\ell_p^n}\right) \subset \mathcal{T}_{w_{0}}^{1,0}\left(\partial B_{\ell_p^n}\right)$;
		\item[(2)] there exists $\lambda\geq 1$ such that
		$$
		{\overline{J_{f}\left(z_{0}\right)}}^{T} v_{w_{0}}=\lambda v_{z_{0}}, \quad \text { or } \quad{\overline{v_{w_{0}}}}^{T} J_{f}\left(z_{0}\right)=\lambda{\overline{v_{z_{0}}}}^{T},
		$$
		where $v_{z_{0}}=\left.\left(\left|z_{1}\right|^{p-2} z_{1}, \dots,\left|z_{n}\right|^{p-2} z_{n}\right)^{T}\right|_{z=z_{0}}$.
	\end{itemize}
\end{theorem}

\begin{remark}
The case $f(z_{0}) = z_{0}$ of Theorem \ref{thm3.1} was studied by Wang and Zhang \cite{WZ24}.
\end{remark}

\begin{remark} Suppose that $f \in H(\mathbb{D}^{n})$ and  $f(0)=0$. In addition, if $f$ is holomorphic at
	$z_{0} \in (\partial \mathbb{D})^{n}$, and satisfies condition
	\[
	f(z_{0})=w_{0}\in (\partial \mathbb{D})^{n},
	\]
	then it is not always true that there exists a constant $\lambda \geq 1$ such that
	\[
	\overline{J_f(z_0)}^{\,T} w_0=\lambda z_0,
	\qquad \text{or equivalently} \qquad
	\overline{w_0}^{\,T}J_f(z_0)=\lambda \overline{z_0}^{\,T}.
	\]
	
	The following example illustrates this fact. Define
	\[
	g(z_1,z_2,\dots,z_{n-1},z_n)
	=
	(z_1^2,z_2,\dots,z_{n-1},z_n)^T.
	\]
	Then $g\in H(\mathbb{D}^{n})$ and $g(0)=0$. Let
	\[
	z_0=w_0=(1,1,\dots,1,1)^T.
	\]
	Clearly,
	\[
	g(z_0)=w_0.
	\]
	
	Further, the Jacobian matrix of $g$ at $z_0$ is given by
	\[
	J_g(z_0)=
	\begin{pmatrix}
		2 & 0 & \cdots & 0 \\
		0 & 1 & \cdots & 0 \\
		\vdots & \vdots & \ddots & \vdots \\
		0 & 0 & \cdots & 1
	\end{pmatrix}.
	\]
	Therefore,
	\[
	J_g(z_0)w_0=
	\begin{pmatrix}
		2 \\
		1 \\
		\vdots \\
		1
	\end{pmatrix},
	\]
	which is not a scalar multiple of $z_0$. Hence, there does not exist any $\lambda \geq 1$ such that
	\[
	\overline{J_g(z_0)}^{\,T}w_0=\lambda z_0.
	\]
	This shows that in the case of polydisc the assertion of Theorem \ref{thm3.1}(2) fails.
\end{remark}
%
We now generalize Theorem C in the setting of pluriharmonic mappings on $B_{\ell_p^n}$ for $p \in [2,\infty]$.
\begin{theorem}\label{thmph}
	Let $p\in[2,\infty]$ and $f: B_{\ell_p^n} \rightarrow B_{\ell_p^N}$ be a pluriharmonic mapping for $n, N \geq 1$. Assume that $f$ is $C^{1}$ at $z_{0} \in \partial B_{\ell_p^n}$ and $f\left(z_{0}\right)=w_{0}$, where $w_{0} \in \partial B_{\ell_p^N}$ and  $w_{0}^{\prime} \in \partial B_p^{2N}$. Let
	$$
	\mathcal{V}_{w_{0}}= \begin{cases}\left(|u_1|^{p-2}u_1,\dots, |u_N|^{p-2}u_N,|v_1|^{p-2}v_1, \dots,|v_N|^{p-2}v_N\right)^T & \text { if } p \geq 2, w_0 \in \partial B_{\ell_p^N} \\
		& \text { and } w_{0}^{\prime} \in \partial B_p^{2N},\\
		 \cfrac{(u,v)^T}{2N} & \text { if } w_0 \in(\partial \mathbb{D})^{N}
	 	{ and } \, w_{0}^{\prime} \in \partial B_{\infty}^{2N},
 \end{cases}
	$$
	where $u=(u_1,\dots,u_N)$ and $v=(v_1,\dots,v_N)$ so that $w_0=(u_1+iv_1,\dots,u_N+iv_N)$. Then
	$$
	\left(J_{f}\left(z_{0}^{\prime}\right)z_{0}^{\prime}\right)^{T} \mathcal{V}_{w_{0}} \geq \frac{1-\left(f(0)^{\prime}\right)^{T} \mathcal{V}_{w_{0}}}{2} \geq \frac{1-\|f(0)^{\prime}\|_p}{2}>0.
	$$
\end{theorem}
Let $X$ be a Banach space and let $B_X$ denote the unit ball in $X$.
Denote by $X^{*}$ the dual space of the real or complex Banach space $X$.
For $x \in X \setminus \{0\}$, define
$$
T(x) = \left\{ \ell_{x} \in X^{*} : \, \ell_{x}(x) = \|x\|_{X} \ \text{ and } \ \|\ell_{x}\|_{X^{*}} = 1 \right\}.
$$
Then, by the well-known Hahn-Banach theorem, we have $T(x) \neq \emptyset$.

\begin{theorem}\label{thm2.7}
	Let $f\in H(\mathbb{D}, B_X)$ and $f(1)\in \partial B_X$. If $f'(1)$ exists, then we have
	\begin{equation*}
		\|f'(1)\|_X \geq \frac{2(1-\|f(0)\|_X)^2}{1-\|f(0)\|_X^2+\|f'(0)\|_X}.
	\end{equation*}
The above estimate is sharp.
\end{theorem}
\begin{remark}
As a consequence of Theorem \ref{thm2.7}, we obtain the sharpness of Theorem E.
\end{remark}
\subsection{ Rigidity Theorem} The following result establishes a rigidity theorem for holomorphic functions on $B_{\ell_p^n}$, $p \in (1,\infty)$.

\begin{theorem} \label{rigi1}
	Suppose that $p \in(1, \infty)$ and that $f\in H(B_{\ell_p^n})$ with $f(0)=0$. For $z \in \partial B_{\ell_p^n}$, let
	$$
	v_{z}=\left(\left|z_{1}\right|^{p-1}, \dots,\left|z_{n}\right|^{p-1}\right)^{T}.
	$$
	If there exist linearly independent real vectors $\alpha_{1}, \alpha_{2}, \dots, \alpha_{n} \in \partial B_{\ell_p^n}$ with non-negative coordinates such that $f$ is holomorphic at $z=\alpha_{k}$ and $f\left(\alpha_{k}\right)=\alpha_{k} \, (k=1, \dots, n)$, then the $n$ equations
	\begin{equation}\label{equal}
	{\overline{v_{\alpha_{k}}}}^{T} J_{f}\left(\alpha_{k}\right) \alpha_{k}=1 \quad(k=1, \dots, n)
	\end{equation}
	hold if and only if $f(z)=z$.
\end{theorem}

\begin{example}
	Consider the mapping
	$$
	f(z)=\bigl(z_{1}z_{n}, \; z_2, \; \dots, \; z_{n-1}, \; z_{n}\bigr)^{T}.
	$$
	We first estimate its norm:
	\begin{align*}
		\|f(z)\|_{p}^{p}
		&= |z_{1}z_{n}|^{p} + |z_{2}|^{p} + \cdots + |z_{n-1}|^{p} + |z_{n}|^{p} \\
		&\leq |z_{1}|^{p} + |z_{2}|^{p} + \cdots + |z_{n-1}|^{p} + |z_{n}|^{p}.
	\end{align*}
	Hence, if \( z \in B_{\ell_p^n} \), then \( \|f(z)\|_{p} \leq 1 \), which shows that \( f \in H(B_{\ell_p^n}) \).
	
	Next, we observe that \( f(0)=0 \) and \( f(e_1)=0 \). Moreover, for each \( j=2,\dots,n \), we have \( f(e_j)=e_j \).
	
	The Jacobian matrix of \( f \) at \( e_j \), for \( j=2,\dots,n-1 \), is given by
	\[
	J_f(e_j) =
	\begin{pmatrix}
		0 & 0 & \cdots & 0 \\
		0 & 1 & \cdots & 0 \\
		\vdots & \vdots & \ddots & \vdots \\
		0 & 0 & \cdots & 1
	\end{pmatrix},
	\]
	while \( J_f(e_n)=I_n \). Consequently,
	\[
	\overline{e_j}^{\,T} J_f(e_j) e_j = 1, \quad \text{for } j=2,\dots,n.
	\]
	
	Despite these properties, \( f(z) \neq z \). This example demonstrates that the assumption in Theorem~\ref{rigi1} requiring \( n \) linearly independent boundary fixed points is essential. In particular, the conclusion \( f(z)=z \) may fail if only \( n-1 \) such fixed points are available.
\end{example}

\begin{example}
	Let \( f(z)=z \) and define \( \alpha_k = -e_k \) for \( 1 \leq k \leq n \). In this case, one can verify that the system of equations in \eqref{equal} is not satisfied.
	
	This example highlights the importance of the assumption that the coordinates of the real vectors \( \alpha_k \), \( 1 \leq k \leq n \), are non-negative.
\end{example}
\section{Proofs of the main results}
 \subsection*{Proof of Theorem \ref{SL}}
	Without loss of generality, we may assume that $\phi(w)=w$ and $z_{0}=(0,\ldots,0)$. Define
	$$
	G(z,w)=\left(
	\frac{f_1(z,w)-w_1}{1-\overline{w_1}f_1(z,w)},\ldots,
	\frac{f_m(z,w)-w_m}{1-\overline{w_m}f_m(z,w)}
	\right).
	$$
	Then $G$ is holomorphic in $z$,  $G(0,w)=0$ and $\|G\|_{\infty}\leq 1$. Consequently, for each $1\leq i\leq m$, the function $G_i(z,w)$ is holomorphic in $z$, vanishes at $z=0$, and $|G_i(z,w)|\leq 1$ for each $i$.
	
	By the Schwarz lemma, we obtain
	$$
	|G_i(z,w)|^{2}\leq \|z\|_{p}^{2}, \qquad 1\leq i\leq m,
	$$
which implies that
	$$
	1-\|z\|_{p}^{2}\leq 1-|G_i(z,w)|^{2}.
	$$
	On the other hand, a direct computation shows that
	$$
	1-|G_i(z,w)|^{2}
	= \frac{(1-|w_i|^{2})(1-|f_i(z,w)|^{2})}{|1-\overline{w_i}f_i(z,w)|^{2}}
	\leq \frac{1-|w_i|^{2}}{|1-\overline{w_i}f_i(z,w)|^{2}}.
	$$
	Combining the above inequalities, we obtain
	$$
	|w_i-f_i(z,w)|^{2}
	\leq |1-\overline{w_i}f_i(z,w)|^{2}
	\leq \frac{1-|w_i|^{2}}{1-\|z\|_{p}^{2}}.
	$$
	By the maximum modulus principle, for $w_i\in r\mathbb{D}$, we have
	$$
	\sup_{w_i\in r\mathbb{D}} |w_i-f_i(z,w)|^{2}
	\leq \frac{1-r^{2}}{1-\|z\|_{p}^{2}}.
	$$
	Letting $r\to 1$ gives that $f(z,w)\equiv w$ and the proof is complete.
\hfill{$\Box$}

\subsection*{Proof of Theorem \ref{thm3.1}}
Let $\alpha \in \mathcal{T}_{z_{0}}(\partial B_{\ell_p^n})$. Without loss of generality, we may assume that $\alpha$ is a unit vector.
Choose a smooth curve $\gamma:[-1,1]\backslash\{0\}\to B_{\ell_p^n}$ such that
$$
\gamma(0)=z_0 \quad \mbox{and} \quad \gamma'(0)=\alpha.
$$
Since $f$ is holomorphic at $z_0$, it follows that $f(\gamma([-1,1])) \subset \overline{B_{\ell_p^n}}$.
Recall that
$$
\rho(z)=\|z\|_p^p-1=|z_1|^p+\cdots+|z_n|^p-1.
$$
Then
$$
\max_{t\in(-1,1)} \rho(f(\gamma(t))) = \rho(f(\gamma(0))) = \rho(w_0)=0.
$$
Hence,
$$
\left.\frac{d}{dt}\rho(f(\gamma(t)))\right|_{t=0}
=2{\rm Re}\left[ \overline{\nabla\rho(w_0)}^{T} J_f(z_0)\alpha \right]=0,
\quad \mbox{ for } \alpha \in \mathcal{T}_{z_{0}}(\partial B_{\ell_p^n}).
$$
This shows that $J_f(z_0)\alpha \in \mathcal{T}_{w_0}(\partial B_{\ell_p^n})$.
Since $J_f(z_0)$ is $\mathbb{C}$-linear, we conclude that
$$
J_f(z_0)\,\mathcal{T}^{1,0}_{z_0}(\partial B_{\ell_p^n}) \subset \mathcal{T}^{1,0}_{w_0}(\partial B_{\ell_p^n}).
$$
Thus, part (1) is proved.

We now show that there exists a $\lambda\in \mathbb{R}$ such that
$$
\overline{J_f(z_0)}^{T} v_{w_0}=\lambda v_{z_0}.
$$
Suppose instead that
$$
\overline{J_f(z_0)}^{T} v_{w_0} = \lambda v_{z_0} + \beta
$$
for some $\lambda\in\mathbb{R}$ and $\beta \in \mathcal{T}_{z_0}(\partial B_{\ell_p^n})$.
From part (1), $J_f(z_0)\beta \in \mathcal{T}_{w_0}(\partial B_{\ell_p^n})$ so that
$$
{\rm Re \,}\langle J_f(z_0)\beta, v_{w_0}\rangle=0.
$$
Since $\beta\in \mathcal{T}_{z_0}(\partial B_{\ell_p^n})$, we also have ${\rm Re \,}\langle \beta, v_{z_0}\rangle=0$.
Therefore,
$$
\|\beta\|_2^2
= {\rm  Re\,}\langle \beta, \lambda v_{z_0} + \beta \rangle
= {\rm Re \,}\langle \beta, \overline{J_f(z_0)}^{T} v_{w_0}\rangle
= {\rm Re \,}\langle J_f(z_0)\beta, v_{w_0}\rangle=0,
$$
which forces $\beta=0$. Thus, $\overline{J_f(z_0)}^{T} v_{w_0}=\lambda v_{z_0}$ with $\lambda\in \mathbb{R}$, and hence
\begin{equation}\label{eq3.1}
	\overline{v_{w_0}}^{T} J_f(z_0)=\lambda \overline{v_{z_0}}^{T}.
\end{equation}

Next, we prove that $\lambda \geq 1$.
Since
\begin{equation}\label{eq3.2}
\frac{\partial \rho(z)}{\partial z}
=\frac{p}{2}\left(|z_1|^{p-2}\overline{z_1},\dots,|z_n|^{p-2}\overline{z_n}\right)
=\frac{p}{2}\,\overline{v_z}^{T},
\end{equation}
the Taylor expansion of $f(z_0-t v_{z_0})$ gives
$$
f(z_0-t v_{z_0})=w_0 - t J_f(z_0)v_{z_0}+O(t^2), \qquad (t\to 0^+).
$$
Using \eqref{eq3.1} and \eqref{eq3.2}, we obtain
\begin{align*}
\|f(z_0-t v_{z_0})\|_p^p -1
&= -p t {\rm Re \,} (\overline{v_{w_0}}^T J_f(z_0)v_{z_0}) + O(t^2)\\
&= -p t \|v_{z_0}\|_2^2 \cdot {\rm Re \,}(\lambda) + O(t^2).
\end{align*}
Since $\lambda\in\mathbb{R}$, we have
\begin{equation}\label{eq3.3}
\|f(z_0-t v_{z_0})\|_p^p = 1 - p\lambda t \|v_{z_0}\|_2^2 + O(t^2).
\end{equation}
On the other hand,
\begin{equation}\label{eq3.4}
\|z_0 - t v_{z_0}\|_p^p = 1 - p t \|v_{z_0}\|_2^2 + O(t^2).
\end{equation}
Finally, Lemma \ref{lem2.3} ensures that
$$
\|f(z_0 - t v_{z_0})\|_p^p \leq \|z_0 - t v_{z_0}\|_p^p, \qquad (t\to 0^+).
$$
Comparing \eqref{eq3.3} and \eqref{eq3.4}, we find that
$$
1 - p\lambda t \|v_{z_0}\|_2^2 + O(t^2) \leq 1 - p t \|v_{z_0}\|_2^2 + O(t^2),
$$
which yields $\lambda \geq 1$.
This completes the proof.
\hfill{$\Box$}
\subsection*{Proof of Theorem \ref{thmph}}
	Let
	\begin{equation}\label{peq3.9}
		\phi(\zeta)=1-\left(f\!\left(\zeta z_{0}\right)^{\prime}\right)^{T} \mathcal{V}_{w_{0}}, \quad \zeta \in \mathbb{D}.
	\end{equation}
Since $\phi(1)=0$, we have $\phi(r)\geq \phi(1)$ for $r \in (0,1)$. Hence,
	$$
	\left(J_{f}\!\left(z_{0}^{\prime}\right) z_{0}^{\prime}\right)^{T} \mathcal{V}_{w_{0}}
	=\lim_{r\to 1^-}\frac{\phi(r)-\phi(1)}{1-r}\geq 0.
	$$
	Since $f$ is pluriharmonic on $B_{\ell_p^n}$, the function $\phi$ defined by \eqref{peq3.9} is non-negative and harmonic on the unit disk $\mathbb{D}\subset\mathbb{C}$.
	By Harnack's inequality in $\mathbb{D}$, it follows that
	$$
	\frac{1-r}{1+r}\,\phi(0)\;\leq\;\phi(\zeta)\;\leq\;\frac{1+r}{1-r}\,\phi(0), \quad |\zeta|=r<1.
	$$
	Consequently,
	$$
	\frac{1}{1+r}\,\phi(0)\;\leq\;\frac{\phi(r)-\phi(1)}{1-r}, \quad 0<r<1.
	$$
	Letting $r \to 1^{-}$ yields
	$$
	\frac{1-\left(f(0)^{\prime}\right)^{T} \mathcal{V}_{w_{0}}}{2}
	\;\leq\;\left(J_{f}\!\left(z_{0}^{\prime}\right) z_{0}^{\prime}\right)^{T} \mathcal{V}_{w_{0}}.
	$$
	 For $ p \in [2,\infty]$, the condition $f(0) \in B_{\ell_p^N}$ implies that $f(0)' \in B_p^{2N}$. Consequently, we obtain
		$$
		\frac{1 - \left(f(0)^{\prime}\right)^{T} \mathcal{V}_{w_0}}{2}
		\;\geq\;
		\frac{1 - \|f(0)^{\prime}\|_p}{2} > 0.
		$$
	Thus, the proof is complete.
\hfill{$\Box$}
\subsection*{Proof of Theorem \ref{thm2.7}} Let $G(\xi) = \ell_b(f(\xi))$, where $\xi \in \mathbb{D}$ and $b = f(1)$.
Then $G \in H(\mathbb{D})$,
$$
G(1) = \ell_b(f(1)) = 1 \quad \mbox{and} \quad G'(1) = \ell_b(f'(1)).
$$
By Theorem D, we have
$$
|G'(1)|\geq \frac{2(1-|G(0)|)^2}{1-|G(0)|^2+|G'(0)|}.
$$
Elementary calculation leads to
\begin{align}\label{eq3.6}
|\ell_b(f(0))|\leq \|f(0)\|_X, \quad |\ell_b(f'(0))|\leq \|f'(0)\|_X\quad \mbox{and} \quad |\ell_b(f'(1))|\leq \|f'(1)\|_X .
\end{align}
Let
$$
g(x)=\frac{2(1-x)}{1+x}=2-\frac{4x}{1+x}, \quad \mbox{for } x\in [0,\infty).
$$
Then, as $ g'(x)=-4/(1+x)^2<0,$ $g$ is a decreasing function in $[0,\infty)$. This fact along with \eqref{eq3.6} gives that
$$
\|f'(1)\|_X\geq \frac{2(1-\|f(0)\|_X)^2}{1-\|f(0)\|_X^2+\|f'(0)\|_X}.
$$
Next, we show the sharpness part. For each fixed $b\in \partial B_X$, we consider
\begin{equation*}
	\phi(z)=b\frac{ A(z)+\|f(0)\|_X}{1+\|f(0)\|_X A(z)},
\end{equation*}
where
$$
A(z)=z \frac{\left(1-\|f(0)\|_X^{2}\right) z+\left\|f^{\prime}(0)\right\|_X}{\left(1-\|f(0)\|_X^{2}\right)+\left\|f^{\prime}(0)\right\|_X z}.
$$
It is easy to see that $\phi \in H(\mathbb{D},B_X)$ and $\phi(1)=b\in \partial B_X$. An elementary calculation shows that
\begin{equation*}
	\|\phi'(1)\|_X = \frac{2(1-\|\phi(0)\|_X)^2}{1-\|\phi(0)\|_X^2+\|\phi'(0)\|_X}.
\end{equation*}
This completes the proof.
\hfill{$\Box$}
\subsection*{Proof of Theorem \ref{rigi1}}
The sufficiency is immediate. It remains to prove the necessity.
Consider
$$
v_{z}=\left(\,|z_{1}|^{p-1}, \ldots, |z_{n}|^{p-1}\,\right)^{T}, \qquad z \in \partial B_{\ell_p^n}.
$$
Then we have
$$
\|v_{\alpha_{k}}\|_{q}=1, \qquad \frac{1}{p}+\frac{1}{q}=1.
$$
Next, we define
$
\phi_{k}(\xi)=v_{\alpha_{k}}^{T} f(\xi \alpha_{k})
$
for $\xi \in \mathbb{D}.$ By H\"older's inequality,
$$
|\phi_{k}(\xi)| \leq \|v_{\alpha_{k}}\|_{q}\,\|f(\xi \alpha_{k})\|_{p} < 1,
$$
which shows that $\phi_{k}\in H(\mathbb{D})$. Moreover, $\phi_{k}(0)=0$ and
$$
\phi_{k}(1)=v_{\alpha_{k}}^{T} f(\alpha_{k})
=v_{\alpha_{k}}^{T}\alpha_{k}
=\|\alpha_{k}\|_{p}^{p}=1.
$$
Furthermore, using \eqref{equal}, we have
$$
\phi_{k}'(1)=v_{\alpha_{k}}^{T} J_{f}(\alpha_{k})\,\alpha_{k}=1.
$$
By Theorem A, we conclude that $\phi_{k}(\xi)=\xi$, i.e.,
$$
v_{\alpha_{k}}^{T} f(\xi \alpha_{k})=\xi.
$$
Differentiating it with respect to $\xi$ and then evaluate it at the origin gives
$$
v_{\alpha_{k}}^{T} J_{f}(0)\,\alpha_{k}=1.
$$
By Lemma \ref{lem2.3}, we have
$$
\|J_{f}(0)\alpha_{k}\|_{p}\leq \|\alpha_{k}\|_{p}=1.
$$
Applying H\"older's inequality again yields
$$
	1 = v_{\alpha_{k}}^{T} J_{f}(0)\alpha_{k}\leq  \|v_{\alpha_{k}}\|_{q} \|J_{f}(0)\alpha_{k}\|_{p}=\|J_{f}(0)\alpha_{k}\|_{p}\leq 1.
$$
Thus, equality holds throughout, and we deduce
$
J_{f}(0)\alpha_{k}=\alpha_{k}.
$
Since the vectors $\alpha_{1}, \ldots, \alpha_{n} \in \partial B_{\ell_p^n}$ are linearly independent, it follows that
$
J_{f}(0)z=z.
$
Finally, by Cartan's rigidity theorem, we obtain
$
f(z)=z.
$
This concludes the result.
\hfill{$\Box$}


\end{document}